\documentclass{amsart}

\input xy
\xyoption{all}

\usepackage[all]{xy}
\usepackage{tikz}
\usepackage{tikz-cd}

\usepackage{enumitem}
\usepackage{geometry}
\usepackage{amsmath}
\usepackage{amssymb}
\usepackage{hyperref}
\usepackage{cite}
\usepackage[final]{showkeys}
\usepackage{hyperref}
\usepackage{setspace}
\usepackage{amsthm}  
\usepackage[pagewise,displaymath,mathlines]{lineno}

\newtheorem{same}{This should never appear}[section]
\newtheorem{defin}[same]{Definition}

\newtheorem{theorem}[same]{Theorem}
\newtheorem{example}[same]{Example}

\newtheorem{fact}[same]{Fact}
\newtheorem{question}[same]{Question}

\newtheorem{prop}[same]{Proposition}

\newtheorem{notation}[same]{Notation}

\newbox\noforkbox \newdimen\forklinewidth
\setbox0\hbox{$\textstyle\smile$}
\setbox1\hbox to \wd0{\hfil\vrule width \forklinewidth depth-2pt
 height 10pt \hfil}
\setbox\noforkbox\hbox{\lower 2pt\box1\lower 2pt\box0\relax}
\def\unionstick{\mathop{\copy\noforkbox}\limits}

\def\nonfork_#1{\unionstick_{\textstyle #1}}

\setbox0\hbox{$\textstyle\smile$}
\setbox1\hbox to \wd0{\hfil{\sl /\/}\hfil}
\setbox2\hbox to \wd0{\hfil\vrule height 10pt depth -2pt width
              \forklinewidth\hfil}
\newbox\doesforkbox
\setbox\doesforkbox\hbox{\lower 2pt\box1 \lower 2pt\box2\lower2pt\box0\relax}
\def\nunionstick{\mathop{\copy\doesforkbox}\limits}

\def\fork_#1{\nunionstick_{\textstyle #1}}

\newcommand{\LS}{\text{LS}}

\newcommand{\ba}{\bold{a}}
\newcommand{\bA}{\bold{A}}

\newcommand{\bK}{\mathbb{K}}
\newcommand{\bx}{\bold{x}}

\newcommand{\rest}{\upharpoonright}

\newcommand{\Mod}{\te{Mod }}

 \newcommand{\footnotei}[1]{}

\newcommand{\te}[1]{\textrm{#1}}

\newcommand{\bL}{\mathbb{L}}

\newcommand{\bY}{\bold{Y}}
\newcommand{\bX}{\bold{X}}
\newcommand{\bZ}{\bold{Z}}

\newcommand{\cC}{\mathcal{C}}

\newcommand{\cF}{\mathcal{F}}

\newcommand{\cL}{\mathcal{L}}

\newcommand{\cP}{\mathcal{P}}

\newcommand{\Set}{\text{Set}}

\renewcommand{\int}{\text{int}}

\newcommand{\ZFC}{\text{ZFC}}
\newcommand{\stat}{\text{stat\, }}

\newcommand{\comment}[1]{}

\title{$\Sigma_1$-Stationary Logic as an $\aleph_1$-Abstract Elementary Class}
\author{Will Boney}
\thanks{\today\\
The author was supported by the National Science Foundation under grants DMS-2137465 and DMS-2339018.\\
We thank Sean Cox for helpful discussions around these ideas and the referee for their thorough report.}

\begin{document}

\begin{abstract}
$\mu$-Abstract Elementary Classes are a model theoretic framework introduced in \cite{bglrv-muaecs} to encompass classes axiomatized by $\bL_{\infty, \infty}$.  We show that the framework extends beyond these logics by showing classes axiomatized in $\bL(aa)$ with just the $aa$ quantifier are an $\aleph_1$-Abstract Elementary Class.
\end{abstract}

\maketitle

\section{Introduction} \label{intro-sec}

In \cite{bglrv-muaecs}, the authors introduced the notion of a $\mu$-Abstract Elementary Class ($\mu$-AECs for short, see Definition \ref{mu-aec-def} below).  $\mu$-AECs are a framework for the model theory of $\bL_{\infty, \mu}$ just as AECs (Abstract Elementary Classes, which are $\aleph_0$-AECs) are a framework for the model theory of $\bL_{\infty, \omega}$. 

A natural question in both cases is to find an example where the relevant logic is weaker than the new framework.  For AECs, the quantifier $Q_{\aleph_1}$ (`there exists uncountably many') and its larger generalizations provides natural examples: there are $\bL(Q_{\aleph_1})$-axiomatizable classes that are not axiomatizable in $\bL_{\infty, \omega}$ but that are AECs when equipped with the natural notion of morphisms\footnote{Examples of this are abound.  The `toy example' of a an equivalence relation with each class countable works, while Kirby \cite[Section 2.8]{k-zpef} shows the more `mathematically interesting' class of Zilber's pseudoexponential structures are another example.}.  Indeed, similar examples exist at higher cardinalities: fixing $\mu<\lambda$, the quantifier $Q_\lambda$ `there exists at least $\lambda$-many' can axiomatize classes not axiomatizable by $\bL_{\infty, \mu}$, but nonetheless these classes form a $\mu$-AEC under the natural choice of morphism (see \cite[Example 8, p. 3051]{bglrv-muaecs}).  However, once `$\mu$' becomes a parameter, these counter examples become less satisfying: $Q_\lambda$ is expressible in $\bL_{\infty, \lambda}$, so a change of parameter (that is already changing) solves the issue.  This motivates the main question of this paper, left implicit in \cite{bglrv-muaecs}:
\begin{question}\label{big-question}
    Is there a $\mu$-Abstract Elementary Class that is not axiomatizable by $\bL_{\infty, \infty}$?
\end{question}

Recent work of Shelah and Villaveces \cite{sh1184} shows that the answer is `no' for $\mu = \aleph_0$. We use stationary logic (introduced by Shelah \cite{sh43}) to give an affirmative answer for $\mu=\aleph_1$ (Theorems \ref{omega1-aec-thm} and \ref{bkm-ex-thm}).  This shows that the framework of $\mu$-Abstract Elementary Classes goes beyond $\bL_{\infty, \infty}$, just as Abstract Elementary Classes go beyond $\bL_{\infty, \omega}$.  We use an example of Barwise, Kauffman, and Makkai \cite[Section 5]{bkm-stationary} to give an explicit class of structures that witness this.

Stationary logic is a fragment of second-order logic that allows for mondaic quantifiers that state that there is a club of realizations in $\cP_{\omega_1}M$ satisfying some formula (or dually that there is a stationary set of realizations), but no other second-order quantifiers; see Section \ref{prelim-ssec} for details.  Importantly, work of Fuchino, Rodrigues, and Sakai \cite{frs-lst-aa-I, frs-lst-aa-II, frs-lst-aa-III} has connected the downward L\"{o}wenheim-Skolem-Tarski behavior of this logic to certain diagonal reflection properties that are independent of $\ZFC$.  Recently, Cox \cite{c-pi1-reflection} has shown that this downward L\"{o}wenheim-Skolem-Tarski property restricted to just instances of the $\stat$quantifier (called there the $\Pi^1_1$ fragment) is already equivalent to such a principle.  In Section \ref{pos-del-sec}, we work with  the $\Sigma_1^1$ fragment (using just positive instances of the $aa$ quantifier) to avoid axioms independent of $\ZFC$.

The proof follows recent results \cite{b-cofquant} on making Abstract Elementary Classes from classes of models axiomatized by the cofinality quantifier.  With that quantifier (as with the `almost all' quantifier $aa$), we wanted to allow only positive instances of the quantifier.  There was an additional issue of `accidental' instances of the quantifier occurring: a countable union of end-extending linear orders will have cofinaly $\omega$ regardless of the cofinality of the indvidual orders.  Thus, we defined a notion of `deliberate' use of quantifiers by deciding for each collection of parameters when want to enforce the quantifier holds for them; this is laid out in Definition \ref{pos-del-def}.

Category theory offers another perspective on Question \ref{big-question}.  The class of $\mu$-AECs over all $\mu$ coincides with the class of accessible categories with monos, \emph{up to equivalence of categories} (see \cite[Theorems 4.3 and 4.10]{bglrv-muaecs}).  In turn, every accessible category is equivalent to a category of models axiomatized by a $\bL_{\infty, \infty}$-theory\footnote{In fact, these theories can be taken to be \emph{basic}, which means each sentence in $T$ is of the form
$$\forall x\left(\phi(\bx)\to \psi(\bx)\right)$$
where $\phi(\bx)$ and $\psi(\bx)$ are positive existential.} by classic results \cite[Theorem 3.2.1 and Corollary 4.3.3]{mp-accessible}.  This seems to preclude an affirmative answer to Question \ref{big-question} if not for the key word `equivalent.' As we previously discussed \cite[Section 4]{bv-structural}, part of this categorical equivalence changes the underlying set of the the models to be unrecognizable from their original presentation presentation.  Responding to a related question of Makkai and Rosicky, Henry \cite{h-nonax-aec} answered a related version of this question by showing that, for each $\mu$, the $\mu$-AEC of `sets of size at least $\mu^+$' is not even equivalent to a $\bL_{\infty, \mu}$-axiomatizable class.  However, it is axiomatizable by $\bL(Q_\mu)$ and, thus, axiomatizable by $\bL_{\infty, \mu^+}$.  For more on this, see Section \ref{pos-del-sec} and Question \ref{inequiv-quest}.

\subsection{Definitions and Notation}\label{prelim-ssec}

We consider several fragments of second-order logic, $\bL^2$.  In this logic, we add second-order variables, quantifiers for these variables\footnote{Although we never consider the second-order universal and existential quantifier in this paper.}, and the atomic formula $`x \in X$' between a first-order variable and a second-order variable.

\begin{notation}
    We use the following conventions:
    \begin{enumerate}
        \item First-order variables are denoted by lowercase letters at the end of the alphabet ($x, y, z, \dots$), and second-order variables are represented similarly by uppercase letters ($X, Y, Z, \dots$).
        \item First-order parameters are denoted by lowercase letters, typically starting at the beginning of the alphabet ($a, b, \dots$), and second-order parameters are denoted by uppercase letters starting at the beginning of the alphabet ($A, B, \dots$).
        \item In all cases, a bold version of the letter is used to denote a string of the appropriate variable length.
        \item In formulas with first- and second-order variables, we use a semicolon to separate them, e.g., $\phi(\bx; \bZ, Y)$ means that $\phi$ is a formula with first-order free variables a subset of the tuple $\bx$ and the second-order free variables are a subset of the tuple $\bZ$  with the singleton $Y$.
    \end{enumerate}
\end{notation}

Following \cite{frs-lst-aa-I}, we define the logic $\bL^{2_{\omega_1}} \subset \bL^2$.  This logic contains first-order logic, allows for no quantification over second-order variables, and restricts the second-order parameters to be countable (or finite) subsets of the universe.  Although this logic has no meaningfully second-order sentences, it is useful to define the desired strong substructure relation in Definition \ref{pos-del-def}.  Crucially, since we restrict to countable sets, this logic is expressible in $\bL_{\omega_1, \omega_1}$, so it inherits many of that logic's desirable qualities (namely, downward L\"{o}wenheim-Skolem-Tarski and closure under $\omega_1$-directed colimits).  We will use (Definition \ref{pos-del-def}.(4)) elementary substructure with respect to $\bL^{2_{\omega_1}}$ as considered by \cite{frs-lst-aa-I}.  Note that \cite[Theorem 1.1.(2)]{frs-lst-aa-I} states that the Downward L\"{o}wenheim-Skolem-Tarski property for $\bL^{2_{\omega_1}}$ (simply called $\cL^{\aleph_0}$ there) to $\aleph_1$ is equivalent to the Continuum Hypothesis since a model could be forced to contain all countable subsets of a countable set.  However, the Downward L\"{o}wenheim-Skolem Tarski property for this logic down to $2^{\aleph_0}$ is a theorem of $\ZFC$.

The main logic we consider is stationary logic, $\bL(aa)$ (under the standard $\omega_1$ interpretation). This is a fragment of second-order logic that interacts with the club filter $\cC(M)$ on $\cP_{\omega_1} M$.  Recall
\begin{eqnarray*}
    \cC(M) = \{X \subset \cP_{\omega_1}M&:& \text{ there is }X_0 \subset X\text{ such that: a) given } A \in \cP_{\omega_1}M\text{ there is }A_0 \in X_0\\
    & &\text{ such that }A \subset A_0\text{; and b) given }A_0 \subset A_1 \subset \dots \text{ from }X_0\text{,}\\
    & &\text{ we have }\cup_{n<\omega}A_n \in X_0\}
\end{eqnarray*}
Alternatively, Kueker \cite[Section 1.1c]{k-count-approx} gives a game-theoretic definitions.  Given a logic $\cL$, $\cL(aa)$ augments it by adding a supply of second-order variables, the atomic formula `$x\in X$', and the quantifier $aa\, s$ over these variables where
\begin{eqnarray*}
    M\vDash `aa \, X \phi(\ba; X, \bA)\text{'} \iff \phi(\ba; M, \bA):=\left\{B \in \cP_{\omega_1}M : M \vDash \phi(\ba; B, \bA)\right\} \in \cC(M)
\end{eqnarray*}

Following \cite[Definition 2.1]{c-pi1-reflection}, we define some fragments of $\cL(aa)$ (which are small parts of a larger hierarchy):
\begin{enumerate}
    \item the \emph{$\Sigma_1$-fragment of stationary logic $\cL(aa)$}, denoted $\cL^{\Sigma_1}(aa)$, consists of all formulas of the form
    $$aa\, Y_1, \dots, aa\,Y_n \psi(\bx; \bX, Y_1, \dots, Y_n)$$
    where $\psi(\bx; \bX, \bY) \in \cL(aa)$ contains no second-order quantifiers; and
    \item the \emph{$\Pi_1$-fragment of stationary logic $\cL(aa)$}, denoted $\cL^{\Pi_1}(aa)$, consists of all formulas of the form
    $$stat\, Y_1, \dots, stat\,Y_n \psi(\bx; \bX, Y_1, \dots, Y_n)$$
    where $\psi(\bx; \bX, \bY) \in \cL(aa)$ contains no second-order quantifiers; and
\end{enumerate}

Any reference to a class-sized logic should be interpreted as referring to some set-sized fragment of it where we do not bother to reference the particular fragment.  For instance, `a theory of $\bL_{\infty, \infty}(\tau)$' really means `a theory of $\bL_{\lambda, \lambda}(\tau)$ for some $\lambda$' (although see Rosicky \cite{r-canon} for some work with proper class sized theories).

We typically use $\cL$ to refer to an abstract logic, and $\bL$ to refer to particular logics.  Thus, `$\bL(aa)$' is the first-order logic $\bL$ augmented by the `almost all' quantifier $aa$, while $`\cL(aa)$' treats $\cL$ as some arbitrary logic (to be plugged in later) that is further augmented by the `almost all' quantifier.

The above describes the `standard $\omega_1$ interpretation' of $aa$, where the club filter is on $\cP_{\omega_1} M$.  We can also vary the $\omega_1$ parameter to define $aa_\lambda$, which restricts second-order variables to subsets of size $<\lambda$ and the semantics of $aa_\lambda$ ask if the definable set is in the club filter on $\cP_\lambda M$.

In our attempts to `Skolemize' $\bL(aa)$, we will make heavy use of the following result of Menas that provides a basis for the various club filters:

\begin{fact}[{\cite{m-strong-super}, see also \cite[Proposition 25.3]{kanamori}}]\label{menas-fact}
Fix $\mu$ and an infinite set $X$.  For each $\cC \subset \cP_\mu X$, $\cC$ contains a club iff there is a function $F:[X]^2 \to \cP_\mu X$ such that
$$\cC(F):=\left\{s \in \cP_\mu X \mid \forall x, y \in s, F(x, y) \subset s \right\} \subset \cC$$
\end{fact}



Finally, we define the main model-theoretic framework we use: $\mu$-Abstract Elementary Classes ($\mu$-AECs).  $\mu$-AECs were introduced in \cite{bglrv-muaecs} as a generalizations of Shelah's Abstract Elementary Classes \cite{sh88}; $\mu$-AECs capture and extend $\bL_{\infty, \mu}$-axiomatizable classes in the same way AECs capture and extend $\bL_{\infty, \omega}$-axiomatizable classes.

\begin{defin}[{\cite[Definition 2.2]{bglrv-muaecs}}] \label{mu-aec-def}
Fix a $<\mu$-ary language $\tau$. Given a class of $\tau$-structures $\bK$ and a binary relation $\prec_\bK$ on $\bK$, we say that $(\bK, \prec_\bK)$ is a \emph{$\mu$-Abstract Elementary Class} iff
\begin{enumerate}
    \item $\prec_\bK$ refines $\tau$-substructure;
    \item $\bK$ and $\prec_\bK$ are both closed under isomorphism;
    \item (Coherence) if $M_0 \subset M_1$ and $M_0 \prec_\bK M_2$ and $M_1 \prec_\bK M_2$, then $M_0 \prec_\bK M_1$;
    \item (Tarski-Vaught Axioms) $\bK$ is closed under $\mu$-directed colimits of $\prec_\bK$ systems and such a colimit is the standard colimit of the structures;
    \item (Downward L\"{o}wensheim-Skolem-Tarski) there is a minimal cardinal $\LS(\bK)$ such that if $M \in \bK$ and $X \subset M$, then there is $N \prec_\bK M$ containing $X$ of size $(X+\LS(\bK))^{<\mu}$.
\end{enumerate}
We refer to $\prec_\bK$ as the `strong substructure relation.'
\end{defin}

$\mu$-Abstract Elementary Classes have a nice characterization in terms of accessible categories where all morphisms are mono.
\begin{fact}[{\cite[Theorems 4.3 and 4.10]{bglrv-muaecs}}]
\begin{enumerate}
    \item Every $\mu$-Abstract Elementary Class is an $(\mu+\LS(\bK))^+$-accesible category where all morphisms are mono.
    \item Every $\mu$-accessible category where all morphisms are monos is equivalent of a $\mu$-Abstract Elementary Class with L\"{o}wenheim-Skolem number $(\mu+\nu)^{<\mu}$, where $\nu$ is the number of morphisms between $\mu$-presentable objects of the category.
\end{enumerate}
\end{fact}

Note that \cite[Definition 3.6]{lrv-internal} gives a more parameterized version of accessibility that better lines up with model-theoretic definitions like $\mu$-AECs

\section{Positive, Deliberate Axiomatizations} \label{pos-del-sec}

Following \cite{b-cofquant}, we give a formal framework for positive, deliberate uses of the $aa$ quantifier.  The framework there was specific to cofinality quantifiers and did not contemplate second order quantifiers, so we must repeat ourselves.

\begin{defin}\label{pos-del-def}
Fix a (first-order) language $\tau$.
\begin{enumerate}
    \item $\tau_* = \tau_*^{\bL^{\Sigma_1}(aa)} = \tau \cup \{R_\phi(\bx; \bZ) \mid \phi(\bx; \bZ, Y) \in \bL^{\Sigma_1}(aa)\}$
    \item The base theory
    $$T^{\bL^{\Sigma_1}(aa)}_{\tau} = \{ \forall \bZ\forall \bx \left( R_\phi(\bx; \bZ) \to aa\, Y \phi(\bx; \bZ, Y) \right) \mid \phi(\bx; \bZ, Y) \in \bL^{\Sigma_1}(aa)\}$$
    \item Given a theory $T \subset \bL^{\Sigma_1}(aa)(\tau)$, define
    \begin{eqnarray*}
        T^* &:=&  \text{``the result of replacing each use of `$aaY\phi(\bx; \bZ, Y)$' with}\\
        & &\text{ `$R_\phi(\bx; \bZ)$' in the inductive construction of each $\psi \in T$"}\\
        T^+ &:=& T^* \cup T^{\bL^{\Sigma_1}(aa)}_\tau
    \end{eqnarray*}
    \item Given $\tau_*$ structures $M \subset N$, we define
    $$M \prec_{(+)} N$$
    to mean that $M$ is an elementary substrcture of $N$ for all formulas in $\bL^{2_{\aleph_1}}(\tau_*)$
\end{enumerate}
\end{defin}

Note that the second-order variables--including the universal quantifier in $T^{\bL^{\Sigma_1}(aa)}_\tau$--are restricted to countable subsets.  Thus, $T^+$ is essentially an $\bL_{\omega_1, \omega_1}(\tau_*)$-theory; this will be important in showing closure under $\omega_1$-directed colimits.

We can now show that we have built an $\omega_1$-AEC:

\begin{theorem}\label{omega1-aec-thm}
Fix $T \subset \bL^{\Sigma_1}(aa)(\tau)$ and fix
$$\bK^+_T =\left( \Mod(T^+), \prec_{(+)}\right)$$
Then $\bK_T^+$ is an $\omega_1$-Abstract Elementary Class with $LS(\bK_T^+) =|T| + \aleph_1$.
\end{theorem}

Note that the $\mu$-AEC axioms only guarantee models of size $\kappa^{<\mu}$ for $\kappa \geq LS(\bK)$, so there is a built-in dependence on the value of the continuum.  

{\bf Proof:} The proof of the various axioms are standard except for the Tarski-Vaught axiom ($\omega_1$-directed colimits) and the downward L\"{o}wenheim-Skolem-Tarski axiom.

\begin{itemize}
    \item {\bf Tarski-Vaught:} Suppose that $\{M_i \in \bK^+_T \mid i \in I\}$ is an $\omega_1$-directed system so $i < j$ implies that $M_i \prec_{(+)} M_j$.  We can build the colimit of this system of $\tau$-structures
    $$M_I := \bigcup_{i \in I} M_i$$
    Since $T^*$ is essentially, $\bL_{\omega_1, \omega_1}$, this all transfers nicely and we are left to show that $M_I$ satisfies $T^{\bL^{\Sigma_1}(aa)}_\tau$.  Fix $\ba \in M_I$; $\bA \in [M_I]^{\aleph_0}$ such that 
    $$M_I \vDash R_\phi(\ba; \bA)$$
    and we want to show that 
    $$\cC:= \{B \in [M_I]^{\aleph_0} \mid M_I \vDash \phi(\ba; \bA, B)\}$$
    contains a club.
    By $\omega_1$-directedness, there is some $i \in I$ so $M_i$ contains the parameters $\ba$ and $\bA$.
    \begin{itemize}
        \item {\bf Unbounded:} Let $B_0 \in [M_I]^{\aleph_0}$.  By $\omega_1$-directedness, find $j > i$ so $B_0 \subset M_j$.  Then $M_j \vDash `R_\phi(\ba; \bA)$' and, thus,
        $$M_j \vDash aa Y \phi(\ba; \bA, Y)$$
        So by unboundedness in $M_j$, there is $B_1 \supset B_0$ in $[M_j]^{\aleph_0}$ such that
        $$M_j \vDash \phi(\ba; \bA, B_1)$$
        Then $B_1 \in [M_I]^{\aleph_0}$ and
        $$M_I \vDash \phi(\ba; \bA, B_1)$$
        as desired.
        \item {\bf Closed:} Let $B_n \in [M_I]^{\aleph_0}$ be increasing so 
        $$M_I \vDash \phi(\ba; \bA, B_n)$$
        Find $i< j \in I$ so $\cup B_n \subset M_j$.  By elementarity, we have for each $n <\omega$
        $$M_j \vDash \phi(\ba; \bA, B_n)$$
        The realizations in $[M_j]^{\aleph_0}$ are club. Then we have that $M_j \vDash `R_\phi(\ba; \bA, \cup B_n)$', so
        $$M_j \vDash \phi(\ba; \bA, \cup B_n)$$
        By $\bL^{2_{\aleph_1}}(\tau_*)$-elementarity, $M_I$ thinks this as well.
    \end{itemize}

    Note that since the structure $M_I$ is computed as the colimit in $\Set$ and $\prec_{(+)}$ is defined nicely, this structure satisfies the second part of the Tarski-Vaught axioms.

    \item {\bf Downward L\"owenheim-Skolem-Tarski:}. Let $M \vDash T^+$.  First, we take advantage of the fact that $T^*$ is essentially $\bL_{\omega_1, \omega_1}$ to (classically) Skolemize it: there is $\tau^{sk}$ and an expansion $M^{sk}$ of $M$ such that, for any $Y \subset M$, we have
    $$\langle Y \rangle_{M^{sk}} \prec_{\bL_{\omega_1, \omega_1}} M$$
    where $\langle \cdot \rangle_{M^{sk}}$ is the closure under the functions of $M^{sk}$.  Notice that the functions of $M^{sk}$ can be $<\omega_1$-ary.
    
    However, we are missing from this conclusion that $\langle Y \rangle_{M^{sk}}$ satisfies $T^{\bL^{\Sigma_1}(aa)}_\tau$; that is, that we have forced the appropriate subsets to contain a club.  For this, we use Menas' result Fact \ref{menas-fact} to further Skolemize to 
    $$\tau^{sk+} := \tau^{sk} \cup \left\{F_\phi:[M]^{2+\ell(\bx) + \omega \cdot \ell(\bZ)} \to [M]^{\aleph_0} \mid \phi(\bx; \bZ, Y) \in \bL^{\Sigma_1}(aa) \right\}$$
    and find a $\tau^{sk+}$-expansion $M^{sk+}$ such that, for all $\ba \in M$ and all $\bA \in [M]^{\aleph_0}$, we have
    $$\left(M\vDash R_\phi(\ba; \bA)\right) \implies \left( \cC\left(F^{M^{sk+}}_\phi(\cdot, \cdot, \ba, \bA)\right) \subset \left\{t \in [M]^{\aleph_0} \mid M \vDash `\phi(\ba; \bA, t)\text{'} \right\}\right)$$
    Above, the interpretation $F_\phi^{M^{sk+}}$ is chosen to witness that $\phi(\ba;\bA, M)$ contains a club.

    Now fix $X \subset M$, and we want to find $M_0 \prec_{(+)} M$ such that $X \subset M_0$ and $\|M_0\| \leq (|T| + \aleph_1 + |X|)^{\aleph_0}$.  Working in the expansion $M^{sk+}$ defined above, set
    $$M_0 := \langle X \rangle_{M^{sk+}}$$
    From this $X \subset M_0$ and is of the desired size.  From above, 
    $$M_0 \prec_{(+)} M$$
    and $M_0 \vDash T^*$.  So it remains to show that $M_0 \vDash T^{\bL^{\Sigma_1}(aa)}_\tau$.

    Fix a formula $\phi(\bx; \bZ, Y)$ and parameters $\ba \in M_0$; $\bA \in [M_0]^{\aleph_0}$, and suppose that 
    $$M_0 \vDash R_\phi(\ba; \bA)$$
    Write $G:[M]^2 \to [M]^{\aleph_0}$ for $F_\phi^{M^{sk+}}(\cdot, \cdot, \ba; \bA)$.  By elementarity, $M \vDash R_\phi(\ba; \bA)$, so this means that 
    $$\cC(G) \subset \left\{B \in [M]^{\aleph_0} \mid M \vDash `\phi(\ba; \bA, B)\text{'} \right\}$$
    Since $M_0$ is closed under $G$, we have that $G \rest [M_0]^2: [M_0]^2 \to [M_0]^{\aleph_0}$.  Let $B \in \cC(G\rest[M_0]^2) \subset [M_0]^{\aleph_0}$.  Then $B \in \cC(G)$, so 
    $$M \vDash \phi(\ba; \bA, B)$$
   By elementarity, this means that $M_0 \vDash \phi(\ba; \bA, B)$.  Thus,
   $$\cC(G\rest [M_0]^2) \subset \{B \in [M_0]^{\aleph_0} \mid M_0 \vDash \phi(\ba; \bA, B)\}$$
   so the latter contains a club, as desired.

   So we have verified each instance of 
   $$M_0 \vDash \forall \bZ \forall \bx\left(R_\phi(\bx;\bZ) \to aa Y \phi(\bx; \bZ, Y)\right)$$
   Thus, $M_0 \in \bK^+_T$, as desired.
\end{itemize}
Putting these together, we have shown that $\bK^+_T$ is an $\omega_1$-AEC.\hfill \dag\\

A version of the Downward LST result appears as \cite[Lemma 4.1]{c-pi1-reflection} by using internal approachability of countable sets.  Our more model-theoretic version using Skolem functions is morally equivalent\footnote{\cite{c-pi1-reflection} has a set-theoretic motivation and structures tend to expansions of some $H_\theta$ (for large regular $\theta$) or appear as constants in such structures.  This means that the countable subsets of $M$--and even club collections of these subsets--also appear as elements of $M$, which changes the argument as subsets also appear as elements.}.  Additionally, our version does not restrict the size of the language.  There is also a more refined Skolemization possible: all second-order variables are quantified by $aa$, so instead of the Skolemization of $aa Y \phi(\bx; \bZ, Y)$ being finitary functions depending on the infinite parameters $\bx; \bZ$, we could incorporate the (future) Skolemization of the $\bZ$ to make the functions entirely finitary.  However, this much more technical formulation provides no advantage in the end result, so we avoid it.\footnotei{WB: But wouldn't it help!  Then we could get models of all sizes?}

We can also generalize this to infinitary logics and other interpretations of $aa$; the proof is a straightforward generalization of the proof of Theorem \ref{omega1-aec-thm} so we omit it.

\begin{theorem}\label{mu-aec-thm}
Fix $T \subset \bL^{\Sigma_1}_{\kappa, \lambda}(aa_\mu)(\tau)$ and set
$$\bK^+_T = \left(Mod(T^+), \prec_{(+_\mu)}\right)$$
Then $\bK^+_T$ is a $(\lambda+\mu)$-Abstract Elementary Class with $LS(\bK^+_T) = |T| + \mu+\kappa$.

In particular, if $T \subset \bL^{\Sigma_1}_{\kappa, \omega_1}(aa)(\tau)$, then $\bK^+_T$ is an $\omega_1$-Abstract Elementary Class.
\end{theorem}

From the set-theoretic work of \cite{frs-lst-aa-I,frs-lst-aa-II,frs-lst-aa-III,c-pi1-reflection}, we know that some on the theories in $\bL(aa)$ that form a $\mu$-AEC is necessary (for a ZFC result at least).  In particular, \cite[Theorem 3.1]{c-pi1-reflection} shows that extending Theorem \ref{mu-aec-thm} to $\bL^{\Pi_1}(aa)$ (which allows for the stationary quantifier) already implies $DRP_{internal}$, a reflection property beyond $\ZFC$.  To prove $DRP_{internal}$, Cox uses this downward L\"{o}wenheim-Skolem property to $\omega_1$ on the following $\Pi_1$-formulas (in the language of set theory):
\begin{itemize}
    \item $$\phi(x) = ``stat\, Z \left( \exists p \left( p = Z \cap \bigcup x \wedge p \in x\right) \right)\text{''}$$
    \item $$\psi(x) = ``stat\, Z \left( \exists p \exists \alpha \left(p = Z \cap x \wedge \alpha < x \wedge \alpha = \sup p\right) \right)\text{''}$$
\end{itemize}

However, the literature is concerned with making the L\"{o}wenheim-Skolem-Tarski property hold down to cardinals below the continuum.  The general question of forming a $\mu$-AEC would allow \emph{any} threshold size, which the current literature doesn't seem to explore.  Thus, the following question is still open:

\begin{question}\label{inequiv-quest}
Is there a theory $T$ from $\bL(aa)$ or $\bL(aa) \cup \bL_{\infty, \infty}$ such that $\Mod(T)$ \emph{cannot} be made into an accessible category for any choice of morphism?

Or is there some set-theoretic hypothesis that implies every $\bL_{\infty, \infty}(aa)$-axiomatizable class can be made into an accessible category via the right choice of morphism?
\end{question}

Sean Cox points out that the existence of a supercompact cardinal $\kappa$ implies a downward L\"{o}wenheim-Skolem-Tarski property for $\bL(aa)$ or even $\bL_{\kappa, \kappa}(aa_\lambda)$ for $\lambda < \kappa$.  However, the closure of these classes under directed colimits is unknown.  In particular, while the reflection of stationary sets is well-studied, the closure of stationary sets under sufficiently directed colimits seems unexplored.

\footnotei{CAN WE SAY ANYTHING INTERESTING ABOUT $\Delta_1$ classes?

\subsection{Equivalence to $\bL_{\infty, \infty}$}}

\section{The Example of Barwise-Kaufmann-Makkai}

To make explicit the claim that there is an $\omega_1$-AEC that is not axiomatizble by $\bL_{\infty, \infty}$, we use an example of Barwise, Kaufmann, and Makkai \cite{bkm-stationary}. Throughout, we will use the following notation: given an ordering $(A, <)$ and subset $X \subset A$, we use
$$(X,<)$$
to denote the obvious ordered set with the inherited linear ordering, perhaps more properly and awkwardly written as $\left(X, <\rest(X\times X)\right)$.

Fix the language
$$\tau = \left(U, V, <, P_U, P_V, E, R\right)$$
where $E, <,$ and $R$ are binary relations and the rest are unary relations.

The key idea of the example is to code two isomorphic well orders $(U,<)$ and $(V,<)$ without exhibiting the actual isomorphism between them.  Instead, $\cP_{\omega_1} U$ and $\cP_{\omega_1} V$ are connected with the $aa$ quantifier to imply that the isomorphism must exist.

First, define $\psi \in\bL_{\omega_1, \omega_1}(\tau)$ that says the following:
\begin{enumerate}
    \item $(U, <)$ and $(V, <)$ are disjoint well orders;
    \item $E$ is an extensional binary relation on $U\times P_U$ and $V \times P_V$;
    \item using $E$, we can identify $P_U$ as $\cP_{\omega_1}U$ and $P_V$ as $\cP_{\omega_1}V$; and
    \item $R$ relates isomorphic pieces of $P_U$ and $P_V$: 
    $$\forall x,y \left[ R(x,y) \leftrightarrow \left(P_U(x) \wedge P_V(y) \wedge \left[(x, <) \cong (y, <)\right]\right)\right]$$
\end{enumerate}

Second, define $\phi\in \bL(aa)$ as
$$aa\, Z \exists x \exists y \left[ P_U(x) \wedge P_V(y) \wedge \left(x =Z \cap U\right) \wedge \left(y = Z \cap V\right) \wedge R(x,y) \right]$$

Two key facts:
\begin{fact}[{\cite[Section 5]{bkm-stationary}}]\label{key-fact}\
\begin{enumerate}
    \item Given $M \vDash \phi$, we have $M \vDash \psi$ iff $(U^M, <^M)$ and $(V^M, <^M)$ are isomorphic.
    \item There is no $\chi \in \bL_{\infty, \infty}(\tau)$ so
    $$\Mod(\chi) = \Mod(\phi \wedge \psi)$$
\end{enumerate}

\end{fact}

\footnotei{[SOMETHING ABOUT THEIR PROOF?]}

Note that Fact \ref{key-fact}.(2) says that $\Mod(\phi\wedge \psi)$ cannot be the set of models of a $\bL_{\kappa, \lambda}$-elementary class for any $\kappa$ and $\lambda$.  This gives us a proof of the following theorem that motivates this article:

\begin{theorem}\label{bkm-ex-thm}
    $\bK_{\phi \wedge \psi}^+$ is an $\omega_1$-Abstract Elementary Class (and thus forms an accessible category) where the models of $\bK^+_{\phi\wedge\psi}$ are not axiomatizable by a $\bL_{\infty, \infty}$-sentence.  Moreover, it is closed under $\omega_1$-directed colimits which are created in $\Set$ and it's $LS$ number is $\omega_1$.
\end{theorem}

{\bf Proof:} From its definition, $\phi \in \bL^{\Sigma_1}(aa)(\tau)$, so $\phi \wedge\psi \in \bL^{\Sigma_1}_{\omega_1, \omega_1}(\tau)$. Then we can use Theorem \ref{mu-aec-thm} for the positive part and Fact \ref{key-fact}.(2) for the negative part.\hfill \dag\\

The abstract theory of accessible category theory implies that $\bK^+_{\phi \wedge \psi}$ is equivalent to a $\bL_{\infty, \omega_1}$-axiomatizable category.  This abstract equivalence is built by `changing the universe', so that a model $M \in \bK^+_{\phi\wedge\psi}$ is represented by the collection of all embeddings from $\omega_1$-sized models of $\bK^+_{\phi\wedge\psi}$ into $M$.  The ability of such a `universe changing' functor to greatly impact the axiomatizability of a class of structures has been observed before (see the discussion in Section \ref{intro-sec}).

However, the Barwise-Kaufmann-Makkai example illustrates something much more subtle.  By being more intentional about turning $\Mod(\phi\wedge\psi)$ into an $\bL_{\infty,\infty}$-axiomatizable class, we can give a very explicit description of the class and an equivalence which shows that the change in axiomatization comes not from changing the universe, but from \emph{a failure of Beth definability}.  This is more akin to the functorial expansions that were used to axiomatize Abstract Elementary classes in \cite[Theorem 3.2.3]{bb-hanf}.

We know that isomorphisms between well-orders are unique.  Thus, the result of \cite{bkm-stationary} implying \emph{existence} of an isomorphism in models of $\phi \wedge \psi$ actually implies the \emph{uniqueness} of the isomorphism; in other words, this isomorphism is implictly definable in $\bL_{\omega_1, \omega_1}(aa)(\tau)$.

On the other hand, Fact \ref{key-fact}.(2) shows that the existence of an isomorphism $U \to V$ is {\bf not} implicitly definable in $\bL_{\infty, \infty}(\tau)$.

We define our strong substructure relation by taking advantadge of this mismatch.  When trying to define a strong substructure $M \prec^+ N$, there are two issues to consider:
\begin{enumerate}
    \item We want it to interact well with $\psi\in \bL_{\omega_1, \omega_1}$.  To do so, let $\cF \subset \bL_{\omega_1, \omega_1}(\tau)$ be the smallest elementary fragment containing the sentence $\phi$; recall this the smallest set of formulas containing $\phi$ and closed under subformulas and all first-order operations.\footnote{In fact, much less than full first-order would suffice, but this is a concise way to capture all of the needed formulas.}  Crucially, if $M \prec_\cF N$, then, for any $a \in P^M_U \cup P^M_V$, we have
    
    $$\{x \in M : M \vDash `xEa\text{'}\} = \{x \in N : N \vDash `xEa\text{'}\}$$
    
    This means that they agree about the isomorphism from $P^M_U$ to $\cP_{\omega_1}(U^M)$ and from $P^M_V$ to $\cP_{\omega_1}(U^V)$.
    \item We want them to agree on the implicit isomorphism defined by $\phi$; see Example \ref{bad-ex} for when this can fail.  To do this, given $M \vDash \phi \wedge \psi$, set $F_M:U^M \to V^M$ to be the unique isomorphism.  We want $M$ being a strong substructure of $N$ to mean that 
    $$F_M = F_N \rest U^N$$
\end{enumerate}

\begin{example}\label{bad-ex}
   We illustrate with a non-example.  Consider the structure $N$ that comes from setting $U^N = V^N = \omega+\omega$ (or use your favorite ordinal in place of $\omega$).  We will build $M \subset N$ by setting $U^M = V^M = \omega$, but we use different parts of $\omega+\omega$: $U^M$ is an initial segment of $U^N$, while $V^M$ is a final segment of $V^N$.  Note that $N \vDash \phi \wedge \psi$ since $U^N \cong V^N$, but that we don't have $M \prec^* N$ since they build incompatible isomorphisms.\\

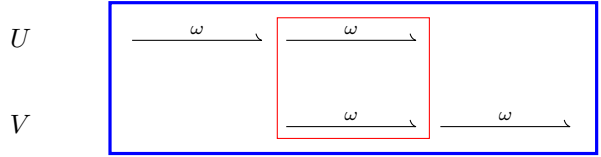
\begin{center}

\begin{tikzcd}[remember picture]
    	U & {} && {} && {} \\
	V &&& {} && {} && {}
	\arrow["\omega", harpoon, from=1-2, to=1-4]
	\arrow["\omega", harpoon, from=1-4, to=1-6]
	\arrow["\omega", harpoon, from=2-4, to=2-6]
	\arrow["\omega", harpoon, from=2-6, to=2-8]
\end{tikzcd}
\end{center}
\begin{tikzpicture}[overlay, remember picture]
    \draw[blue, very thick] (4.9,0) rectangle (11.3,2);
    \draw[red, thin] (7.1,0.2) rectangle (9.1,1.8);
\end{tikzpicture}

The outer box represents $U^N \cup V^N$ and the inner box represents $U^M \cup V^M$.

Thus, $M\prec_\cF N$ but $M \not \prec^+ N$.

\end{example}

To streamline our definition, we define the auxiliary class of models where the implicitly definable isomorphism is explicit.

\begin{defin} \
\begin{enumerate}
    \item Set $\tau^+ = \tau \cup \{F\}$ for a unary function symbol and $\psi_F$ to be the first-order statement that $F: U \to V$ is a bijection.
    \item Set $\bK^+ = \left(\Mod (\phi \wedge \psi_F), \prec_\cF\right)$ and define the map
    \begin{eqnarray*}
    \Mod (\phi \wedge \psi) &\to& \Mod (\phi \wedge \psi_F)\\
    M &\mapsto& M^+ = (M, F_M)
    \end{eqnarray*}
    where $F_M$ is the unique isomorphism $U^M \to V^M$.
\end{enumerate}
\end{defin}
\begin{defin}
Define the class $\bK = \left( \Mod(\phi\wedge\psi), \prec^+\right)$, where $\prec^+$ is defined by, for any $M, N \vDash \phi \wedge\psi$,
\begin{eqnarray*}
    M \prec^+ N & \iff& M^+ \prec_\cF N^+
\end{eqnarray*}
    
\end{defin}

We now collect several straightforward results.
\begin{prop} \
\begin{enumerate}
    \item $\bK^+$ is $\bL_{\omega_1, \omega_1}$-axiomatizable and, therefore, an $\omega_1$-Abstract Elementary Class.
    \item The map $M \mapsto M^+$ induces an isomrophism (and thus an equivalence) between $\bK$ and $\bK^+$ that preserves the underlying sets (a \emph{concrete} isomorphism).
    \item $\bK$ is an $\omega_1$-Abstract Elementary Class, but models of $\bK$ are {\bf not} axiomatizable by any $\bL_{\infty, \infty}$ sentence.
\end{enumerate}    
\end{prop}

{\bf Proof:} (1) follows because $\cF$ is a fragment of $\bL_{\omega_1, \omega_1}$ containing the relevant sentences (see item (3) in \cite[Section 2]{bglrv-muaecs}).  (2) is the content of saying the bijection $F_M$ is implicitly definable from $\phi \wedge \psi$. (3) follows from (1) and (2).  \hfill \dag\\

Work in progress with Mick Walker reveals that there is a deep connection between interpolation (seen as a property that implies Beth definability) and the axiomatization of $\mu$-Abstract Elementary Classes.

\bibliography{bib}
\bibliographystyle{alpha}

\end{document}